

\baselineskip=14pt
\parskip=10pt

\magnification=\magstephalf

\def\1{{\overline{1}}}
\def\2{{\overline{2}}}
\parindent=0pt
\overfullrule=0in

\def\frac#1#2{{#1 \over #2}}
\centerline
{\bf 
Implementing and Experimenting with the Calabi-Wilf algorithm 
}
\centerline
{\bf for random selection of a subspace over a finite field}
\bigskip
\centerline
{\it Shalosh B. EKHAD and Doron ZEILBERGER}
\bigskip
\qquad \qquad \qquad {\it  In fond memory of Eugenio Calabi (May 11, 1923 - September 25, 2023), and Herbert Saul Wilf (June 13, 1931 - January 7, 2012)}

{\bf Abstract:} 
In their beautiful note, ``On the Sequential and Random Selection of Subspaces over a Finite Field", geometrical giant Eugenio Calabi and combinatorial giant Herbert 
Wilf proposed an elegant algorithm to do what is promised in their title. In the present note, written forty six years later, 
we describe a Maple package, written by the second author, that implements their beautiful algorithm, and then go on to report numerous experiments, performed by the first author, that demonstrate the efficiency and reliability of the Calabi-Wilf algorithm.

{\bf Memory Lane (by DZ)}

One of my numerous temporary positions was during the academic year 1982-1983. Herbert Wilf, who later became my close collaborator, got me
a one-year visiting position at the University of Pennsylvania. I occupied the office of David Harbater, who was on sabbatical that year, and it so happened
that it was  close to {\bf both} Wilf's and Calabi's offices. I was in awe of both of them,  since they were already legendary way back then.
One day Herb mentioned to me that, while Calabi's expertise was almost {\it diametrically opposite} to his,
they did have a {\it non-empty intersection}, and coauthored a cute little note [CW]  that applies the  {\it Wilf methodology} [W]
for random selection of {\it combinatorial} objects, to something that besides being  combinatorial (after all `finite fields' are finite!) was also
a {\it geometrical} object, that Calabi can relate to: namely {\bf subspaces} of $GF(q)^n$.

{\bf Reminder about the Wilf approach for Random Generation of Combinatorial objects}

The {\it binomial coefficient},
$$
{{n} \choose {k}} = \frac{n!}{k!(n-k)!} \quad,
$$
famously counts the number of {\it subsets} with $k$ elements of a `universal' set of $n$ elements,$\{1, \dots, n\}$. It famously satisfies the
{\it Pascal-Chu} recurrence
$$
{{n} \choose {k}} = {{n-1} \choose {k-1}} \, + \,   {{n-1} \choose {k}} \quad,
$$
that {\it breaks up} all $k$-element subsets of $\{1, \dots, n\}$ into 
those that {\bf do} contain $n$, whose number is the first term on the right, and those that
{\bf do not} contain $n$, whose number is the second term on the right.

It follows that the {\it fraction} of $k$-element subsets of $\{1,\dots,n\}$ that {\bf do} contain $n$ is
$$
\frac
{{{n-1} \choose {k-1}}}{{{n} \choose {k}}}
= 
\frac{(n-1)!/((k-1)!(n-k)!)}{(n!/k!(n-k)!)}
=
\frac{k}{n} \quad ,
$$
hence the {\it probability} that a random $k$-subset of $\{1, \dots, n\}$ {\bf does} contain $n$ is $\frac{k}{n}$, while the probability that
it {\bf does not} is $1- \frac{k}{n}$.

This simple observation lead Wilf [W] to suggest the following {\it recursive} algorithm for the {\it random generation},  {\bf uniformly at random}, of a $k$-subset of $\{1,\dots, n\}$.

{\tt RandomSubset(n,k)}: if $k<0$ or $k>n$ return {\tt FAIL}. If $k=0$ then RETURN the empty set. Otherwise
toss a loaded coin whose probability of Heads is $k/n$. If it lends on Heads, output
{\tt RandomSubset(n-1,k-1)}$\cup \{n\}$, else output {\tt RandomSubset(n-1,k)} .

This {\it philosophy} was used on many other combinatorial objects that posses {\it recursive structures}, and nicely described, and {\it implemented} (with FORTRAN source-code included!), in the
Nijenhuis-Wilf classic [NW]. However, the Calabi-Wilf algorithm did not make it to that book. It was exposited in more detail by Igor Pak,
in his insightful `big picture'  paper [P], written in his inimitable engaging style.

{\bf The Calabi-Wilf algorithm}

The {\it q-binomial coefficient},
$$
\left [ \matrix{n \cr k} \right ]=\frac{[n]!}{[k]![n-k]!} \quad,
$$
(where $[n]! := (1-q) \cdots (1-q^n)/(1-q)^n$ is the $q$-analog of $n!$),
not-quite-as-famously, counts the number of $k$-dimensional {\it subspaces} of  $GF(q)^n$.
Each such subspace has many possible bases, of course, but only one in {\it row-echelon} form, so in order to count $k$-dimensional subspaces
one should count $k \times n$ matrices over $GF(q)$ in row-echelon form.

The right side of the {\it q-Pascal-Chu} recurrence
$$
\left [ \matrix{n \cr k} \right ]\,=\,
\left [ \matrix{n-1 \cr k-1} \right ] \, + \, q^k \left [ \matrix{n-1 \cr k} \right ]  \quad,
$$
may be thought of (per [CW]) as counting complementary subsets of the totality of $k \times n$ matrices in echelon form:

$\bullet$ The first term counts those $k \times n$ matrices,  $B=(b_{ij})_{1 \leq i \leq k, 1 \leq j \leq n}$ for which $b_{11}=1$, $b_{i1}=0$ ($1<i \leq k$), $b_{1j}=0$ ($1<j \leq n)$, and
for which {\it necessarily} the remaining $(k-1) \times (n-1)$ array is a $(k-1) \times (n-1)$ matrix in row-echelon form.

$\bullet$ The second term counts the other basis matrices $B$, i.e. those whose first column is an {\bf arbitrary} vector of length $k$, and for which the
remaining $k \times (n-1)$ array is a basis matrix for a $k$ dimensional subspace of $GF(q)^{n-1}$.

It follows that the {\it fraction} of matrices of the first kind is
$$
\left [ \matrix{n-1 \cr k-1} \right ]/\left [ \matrix{n \cr k} \right ]= \frac{q^k-1}{q^n-1} \quad,
$$
and it follows, in turn, that the {\it probability} that a random $k\times n$ matrix in row-echelon form is of the first kind is
$\frac{q^k-1}{q^n-1}$.

This  observation lead Calabi and Wilf [CW] to suggest the following {\it recursive} algorithm for the {\it random generation},  {\bf uniformly at random}, of a $k \times n$ matrix over $GF(q)$ in row-echelon form.

{\tt RandomSubspace(n,k)}: if $k<0$ or $k>n$ return {\tt FAIL}. If $k=0$ then RETURN the empty matrix (corresponding to the null subspace). Otherwise
toss a loaded coin whose probability of Heads is $(q^k-1)/(q^n-1)$. If it lends on Heads, generate
{\tt RandomSubspace(n-1,k-1)}, and convert it to a $k \times n$ matrix by putting on top the row vector $10^{n-1}$, and  prefixing the column vector
$(10^{k-1})^T$ (where of course the sole $1$ is shared by both).

Else generate {\tt RandomSubspace(n-1,k)}. Also generate a random column vector in $GF(q)$, of length $k$, and prefix it at the front.

{\bf Maple Implementation}

The Calabi-Wilf algorithm (and much more!) is implemented in the Maple package, \hfill\break
{\tt CalabiWilf.txt}, accompanying the present paper:

{\tt https://sites.math.rutgers.edu/\~{}zeilberg/tokhniot/CalabiWilf.txt}  \quad.

It is procedure {\tt Rqnk(q,n,k)}. For example, typing

{\tt Rqnk(7,10,5);}

would generate, {\it uniformly-at-random}, one of the 
$$
((7^{10}-1)(7^{9}-1)\cdots (7^6-1))/((7^5-1)(7^4-1) \cdots (7-1))=1602592475815614015216
$$ 
$5 \times 10$ row-echelon matrices over $GF(7)$. Here is
one of them:
$$
\left (
\matrix{
5 & 6 & 3 & 2 & 5 & 1 & 0 & 0 & 0 & 0 \cr
 6 & 1 & 0 & 0 & 6 & 0 & 1 & 0 & 0 & 0 \cr
 5 & 3 & 4 & 2 & 0 & 0 & 0 & 1 & 0 & 0 \cr
 0 & 5 & 1 & 6 & 6 & 0 & 0 & 0 & 1 & 0 \cr
 5 & 6 & 2 & 1 & 5 & 0 & 0 & 0 & 0 & 1 
}
\right ) \quad .
$$

{\bf Experimentation}

$\bullet$ If you want to see an article about estimating the average, variance, skewness,  and kurtosis,
for the number of occurrences of $1$s, in row-echelon matrices of dimension $k \times 2k$ 
over GF(2), for $k$ from $50$ to  $100$, by simulating, $1000 $ times, 
look at the output file

{\tt http://www.math.rutgers.edu/\~{}zeilberg/tokhniot/oCalabiWilf1.txt} \quad .

$\bullet$ If you want to see an article about estimating the average, variance, skewness,  and kurtosis, for the number of occurrences (as consecutive submatrix) of the matrix
[[1, 0, 2], [1, 0, 2], [1, 0, 1]],  in row-echelon matrices of dimension $k \times 2k$  over $GF(3)$ for $k$ from $50$ to $60$ by simulating $1000$ times, look at the output file

{\tt http://www.math.rutgers.edu/\~{}zeilberg/tokhniot/oCalabiWilf1a.txt} \quad .

$\bullet$  If you want to see an article about estimating  the average, variance, skewness, and kurtosis
of the number of $1$s in $k \times 2k$ row-echelon matrices over $GF(3)$ for
$50\leq k\leq 55$ , by simulating $1000$ times (and repeating each run three times, to compare notes), and comparing the estimated values to the (theoretically computed) exact values, look at

{\tt http://www.math.rutgers.edu/\~{}zeilberg/tokhniot/oCalabiWilf2.txt} \quad .

$\bullet$  If you want to see an article about estimating  the average
of the number of $1$s in $k \times 3k$ row-echelon matrices over $GF(2)$ for
$100\leq k\leq 110$ , by simulating $1000$ times (and repeating each run three times, to compare notes), and comparing the estimated average to the (theoretically computed) exact value, look at

{\tt http://www.math.rutgers.edu/\~{}zeilberg/tokhniot/oCalabiWilf3.txt} \quad .

$\bullet$ If you want to see an article on the average minimal weight of vectors in $k$-dimensional subspaces of $GF(2)^n$ for $k$ from $1$ to, $5$, and $n$ from  $10$ to  $100$, in increments of $10$,
look at the output file

{\tt http://www.math.rutgers.edu/\~{}zeilberg/tokhniot/oCalabiWilf4.txt} \quad .

For more output files see the web-page of this paper:

{\tt https://sites.math.rutgers.edu/\~{}zeilberg/mamarim/mamarimhtml/calabi.html} \quad .

Readers are welcome to experiment with our Maple package to their heart's content, and generate many more output files.

{\bf Conclusion}: Rest in peace, my dear heroes, Herb and Gene, you both did great work on your own, and in collaboration with other people, in very different parts
of mathematics, but this sole joint work is also great!

{\bf References}

[CW] Eugenio Calabi and Herbert S. Wilf, {\it On the Sequential and Random Selection of Subspaces over a Finite Field}, 
J. Of Combinatorial Theory (Ser. A) {\bf 22} (1977), 107-109. \hfill\break
Reprinted in:
Eugenio Calabi, {\it ``Collected works''}, Edited by Jean-Pierre Bourguignon, Xiuxiong Chen and Simon Donaldson. With contributions by Shing-Tung Yau, Blaine Lawson, Marcel Berger and Claude LeBrun
Springer, Berlin, (2020), 481-483.
\hfill\break
{\tt https://sites.math.rutgers.edu/\~{}zeilberg/akherim/CalabiWilf1977.pdf}  \quad .

\vfill\eject

[NW] Albert Nijenhuis and Herbert S. Wilf, {\it ``Combinatorial Algorithms''}, Academic Press. First edition: 1975. Second edition: 1978.  \hfill\break
{\tt https://www2.math.upenn.edu/\~{}wilf/website/CombinatorialAlgorithms.pdf} \quad .

[P] Igor Pak, {\it When and how n choose k},  {\it Randomization methods in algorithm design}, 191-238, 1997.  \hfill\break
{\tt https://www.math.ucla.edu/\~{}pak/papers/nk13.pdf} \quad.

[W] Herbert S. Wilf, {\it A unified setting for sequencing, ranking, and selection algorithms for combinatorial objects}, 
Advances in Math {\bf 24} (1977) , 281-291. \hfill\break
{\tt https://www2.math.upenn.edu/\~{}wilf/website/Unified\%20setting.pdf} \quad .

\bigskip
\hrule
\bigskip

Shalosh B. Ekhad, c/o D. Zeilberger, Department of Mathematics, Rutgers University (New Brunswick), Hill Center-Busch Campus, 110 Frelinghuysen
Rd., Piscataway, NJ 08854-8019, USA. \hfill\break
Email: {\tt ShaloshBEkhad at gmail  dot com}   \quad .
\smallskip

Doron Zeilberger, Department of Mathematics, Rutgers University (New Brunswick), Hill Center-Busch Campus, 110 Frelinghuysen
Rd., Piscataway, NJ 08854-8019, USA. \hfill\break
Email: {\tt DoronZeil at gmail  dot com}   \quad .
\bigskip
Oct. 27, 2023.

\end